\documentclass{amsart}

\newcommand{\bbf}{\mathbb{F}}

\newcommand{\bbp}{\mathbb{P}}\newcommand{\bbr}{\mathbb{R}}\newcommand{\bbz}{\mathbb{Z}}
\newcommand{\caa}{\mathcal{A}}
\newcommand{\cae}{\mathcal{E}}

\newcommand{\cas}{\mathcal{S}}

\newcommand{\frD}{\mathfrak{D}}

\newcommand{\ch}{\hat{h}}

\newcommand{\ov}{\overline}

\DeclareMathOperator{\tor}{tor}

    \newtheorem{theorem}{Theorem}
\numberwithin{theorem}{section} \theoremstyle{plain}

\newtheorem{corollary}[theorem]{Corollary}

\newtheorem{lemma}[theorem]{Lemma}

\newtheorem{proposition}[theorem]{Proposition}

\theoremstyle{definition}

\newtheorem{remark}[theorem]{Remark}
\numberwithin{equation}{section}

\begin{document}
\title{Analogues of Lehmer's problem in positive characteristic}

\author{Am\'{\i}lcar Pacheco}

\date{\today}

\address{Universidade Federal do Rio de Janeiro (Universidade do Brasil)\\
Departamento de Ma\-te\-m\'a\-ti\-ca Pura\\
Rua Guai\-aquil 83, Cachambi, 20785-050 Rio de Janeiro, RJ, Brasil}

\thanks{This work was partially supported by CNPq research grant 300896/91-3 and Pronex \#41.96.0830.00.}

\email{amilcar@impa.br}


\begin{abstract}Let $C$ be a smooth projective irreducible curve defined over a finite field
$\bbf_q$ and $K=\bbf_q(C)$. We show that every non-torsion element $\alpha\in\ov{K}$ of degree $d$ over $K$ of a
Drinfeld $A$-module $\phi$ defined over $K$ has canonical height $\ch_{\phi}(\alpha)$ at least $1/d$. Similarly, if
$E/K$ is a non-constant elliptic curve defined over a function field $K=l(C)$ of a curve $C$ defined over an
algebraically closed field $l$ of characteristic 0 or $p>3$, we show that every point of infinite order $P\in
E(\ov{K})$ of degree $d$ over $K$ has canonical height $\ch_E(P)$ at least $c/d$, where $c$ depends only on the degree
of the $j$-map associated to $E/K$.
\end{abstract}

\maketitle

\section{Introduction}
Let $\alpha$ be an algebraic number of degree $d$ over $\mathbb{Q}$ and suppose it is not a root of unity. Let
$h:\overline{\mathbb{Q}}\to\mathbb{R}$ be the absolute logarithmic height. Lehmer's conjecture consists in asking for
an absolute real constant $c>0$ such that $h(\alpha)\ge\frac cd$.

Although this question remains open, analogues of this conjecture have been considered in other contexts. Let $E$ be an
elliptic curve defined over a number field $K$, $j_E$ its $j$-invariant, $\overline{K}$ the algebraic closure of $K$,
$P\in E(\overline{K})$ a point of infinite order and $\hat{h}_E:E(\overline{K})\to\mathbb{R}$ its canonical height. Let
$K(P)$ be the field generated over $K$ by the coordinates of $P$, $d=[K(P):K]$ and $D=[K:\mathbb{Q}]$. In
\cite[Corollary 0.2]{HiSi} it is shown that if $j_E$ is non-integral, then there exists $c>0$ depending on $E/K$ such
that $\hat{h}_E(P)\ge\frac c{d^2(\log d)^2}$. Let $h=\max\{1,h(j_E)\}$. This result was improved in \cite[Corollary
1.4]{Da}, where it was proved that there exist absolute effective computable real constants $c_5,c_6>0$ such that
$\hat{h}_E(P)\ge c_5h(dD)^{-3}\left(1+\frac{\log(dD)}h\right)^{-2}$, if $j_E$ is integral, and $\hat{h}_E(P)\ge
c_6D^{-3}d^{-15/8}h^{-2}\left(1+\frac{\log(dD)}h\right)^{-2}$, otherwise. In section 3 we prove an analogue of this
result for non-constant elliptic curves over function fields over algebraically closed fields of characteristic 0 or
$p>3$.

Let $C$ be a smooth irreducible projective curve defined over a finite field $\mathbb{F}_q$ of $q$ elements and
$K=\mathbb{F}_q(C)$. The direct translation of Lehmer's conjecture to $K$ is trivial, because the requirement that
$\alpha$ is not a root of unity is equivalent to $\alpha\in K-\mathbb{F}_q$, thus there is a discrete valuation
$v:K\to\mathbb{Z}\cup\{\infty\}$ of $K$ such that $v(\alpha)<0$ and therefore $h(\alpha)\ge\frac 1d$.

Another instance of the Lehmer problem is to consider the canonical height $\ch_{\phi}:\ov{K}\to\bbr$ of a Drinfeld
$A$-module $\phi:A\to K\{\tau\}$ of rank $r$ defined over a $K=\bbf_q(C)$. We also have a notion of torsion elements in
this context and we ask for a constant $c$ depending on $\phi$ such that for every non-torsion element
$\alpha\in\ov{K}$ with $[K(\alpha):K]=d$ we have $\ch_{\phi}(\alpha)\ge\frac cd$. We prove this in section 2 starting
with the case where $K=\bbf_q(T)$ is the rational function field over $\bbf_q$ and $A=\bbf_q[T]$. The general result is
then deduced from this case. Analogues of this type of result were proved in \cite{De1}.

\section{Drinfeld modules}
Let $A=\bbf_q[T]$ be the polynomial ring in one variable over the finite field $\bbf_q$ of $q$ elements, $k=\bbf_q(T)$
its field of fractions and $\phi:A\to\text{End}_k(\mathbb{G}_a)\cong k\{\tau\}$ a Drinfeld $A$-module of rank $r$
defined over $k$ with respect to the inclusion $A\subset k$. Denote $\phi_T=T+a_1\tau+\ldots+a_r\tau^r$.

Let $\ov{k}$ be the algebraic closure of $k$ and $h:\ov{k}\to\bbr$ the absolute logarithmic Weil height. The global
height of the Drinfeld module $\phi$ at $\alpha\in\ov{k}$ is defined by (cf. \cite[\S2]{De})
$$\ch_{\phi}(\alpha)=\lim_{n\to\infty}\frac{h(\phi_{T^n}(\alpha))}{q^{nr}}.$$
Analogously to the case of elliptic curves this global height decomposes in a sum of local heights which are defined as
follows. Let $L=k(\alpha)$ and $M_L$ the set of places of $L$ normalized so that they correspond to discrete valuations
$v:L\to\mathbb{Z}\cup\{\infty\}$. Let $d_v$ be the degree of $v$ and $d=[L:k]$. The local height of $\alpha$ at $v$
with respect to $\phi$ is defined by (cf. \cite[\S4]{Po})
$$\ch_{\phi,v}(\alpha)=-\frac{d_v}d\lim_{n\to\infty}\frac{\min\{0,v(\phi_{T^n}(\alpha))\}}{q^{nr}}.$$

It follows from the above definitions that
\begin{equation}\label{ht3}
\ch_{\phi}(\alpha)=\sum_{v\in M_L}\ch_{\phi,v}(\alpha).
\end{equation}
An element $\alpha\in\ov{k}$ is called a torsion element of $\phi$ if there exists $f\in A-\{0\}$ such that
$\phi_f(\alpha)=0$.

\begin{theorem}\label{thmdm}
Let $\alpha\in\ov{k}$ be a non-torsion element of $\phi$ and $d=[k(\alpha):k]$. Then
$$\ch_{\phi}(\alpha)\ge\frac1d.$$
\end{theorem}

\begin{proof}
Let $S\subset M_L$ be the set consisting of the poles of $T=a_0$, $a_1,\cdots,a_r$ and the zeros of $a_r$. Suppose
there exists $v\notin S$ such that $v(\alpha)<0$. Then, for every $0\le i<d$,
$$q^rv(\alpha)+v(a_r)=q^rv(\alpha)<q^iv(\alpha)\le q^iv(\alpha)+v(a_i),$$
hence $v(\phi_T(\alpha))=q^rv(\alpha)$. By induction, for every $n\ge 1$, we also have
$v(\phi_{T^n}(\alpha))=q^{nr}v(\alpha)$, thus $\ch_{\phi,v}(\alpha)=-\frac{d_v}dv(\alpha)\ge\frac1d$.

Assume now that all the poles of $\alpha$ lie in $S$. Let $v\in S$ be a pole of $\alpha$. Let
$$M_{\phi,v}=\min_{0\le i<r}\frac{v(a_i)-v(a_r)}{q^r-q^i}.$$
Suppose $v(\alpha)<M_{\phi,v}$ and $v(a_r)\le 0$. The first inequality implies $v(\phi_T(\alpha))=v(a_r)+q^rv(\alpha)$.
The two inequalities imply
$$q^r(q^r-q^i)v(\alpha)<(q^r-q^i)v(\alpha)<v(a_i)-v(a_r)\le v(a_i)-v(a_r)-(q^r-q^i)v(a_r),$$
for every $0\le i<r$, i.e.,
$$q^{2r}v(\alpha)+(q^r+1)v(a_r)<q^{r+i}v(\alpha)+q^iv(a_r)+v(a_i),$$
i.e.,
$$v(\phi_{T^2}(\alpha))=q^{2r}v(\alpha)+(q^r+1)v(a_r)=q^{2r}v(\alpha)+\frac{q^{2r}-1}{q^r-1}v(a_r).$$
Suppose we have proved that for every integer $1\le m<n$ we have
$$v(\phi_{T^m}(\alpha))=q^{mr}v(\alpha)+\frac{q^{mr}-1}{q^{r}-1}v(a_{r}).$$
Then
\begin{equation*}\begin{aligned}
&q^{r}(q^{r}-q^i)v(\phi_{T^{n-2}}(\alpha))\\
&=q^{r}(q^{r}-q^i) (q^{(n-2)r}v(\alpha)+(q^{(n-3)r}+\ldots+q^{r}+1)v(a_{r}))\\
&\le q^{(n-1)r}(q^{r}-q^i)v(\alpha)<(q^{r}-q^i)v(\alpha)<v(a_i)-v(a_{r})\\
&\le v(a_i)-v(a_{r})-(q^{r}-q^i)v(a_{r}).
\end{aligned}\end{equation*}
Thus,
\begin{equation*}\begin{aligned}
&q^{nr}v(\alpha)+(q^{(n-1)r}+\ldots+q^{r}+1)v(a_{r})\\
&<q^{(n-1)r+i}v(\alpha)+(q^{(n-2)r+i}+\ldots+q^{r+i}+q^i)v(a_{r})+v(a_i),
\end{aligned}\end{equation*}
i.e.,
$$v(\phi_{T^n}(\alpha))=q^{nr}v(\alpha)+(q^{(n-1)r}+\ldots+q^{r}+1)v(a_{r})=
q^{nr}v(\alpha)+\frac{q^{nr}-1}{q^{r}-1}v(a_{r}).$$ Hence,
$$\ch_{\phi,v}(\alpha)=-\frac{d_v}d\left(v(\alpha)+\frac1{q^r-1}v(a_r)\right)\ge\frac1d.$$

Suppose now that $v(\alpha)<M_{\phi,v}$, but $v(a_r)>0$. Let $\xi$ be a sufficiently negative power of a local
parameter at $v$ so that $v(\xi^{q^r-1}a_r)\le0$. The Drinfeld module $\psi=\xi^{-1}\phi\xi$ is isomorphic to $\phi$
and by \cite[Proposition 2]{Po} $\ch_{\phi,v}=\ch_{\psi,v}$. Note that
$\psi_T=T+\xi^{q-1}a_1\tau+\ldots+\xi^{q^r-1}a_r\tau^r$. Then for every $0\le i<r$ we have
$$(q^r-q^i)v(\alpha)<v(a_i)-v(a_r)<v(a_i)-v(a_r)-v(\xi)(q^r-q^i)=v(\xi^{q^i-1}a_i)-v(\xi^{q^r-1}a_r),$$
in particular, $v(\alpha)<M_{\psi,v}$. By the argument of the last paragraph we conclude that
$\ch_{\phi,v}(\alpha)=\ch_{\psi,v}(\alpha)\ge\frac1d$.

If $v(\alpha)\ge M_{\phi,v}$, let $\xi$ be a sufficiently positive power of a local parameter at $v$ such that
$$M_{\psi,v}=\min_{0\le i<r}\frac{v(a_i\xi^{1-q^i})-v(a_r\xi^{1-q^r})}{q^r-q^i}
=\min_{0\le i<r}\left(\frac{v(a_i)-v(a_r)}{q^r-q^i}+v(\xi)\right)>v(\alpha).$$ Once again we take the Drinfeld module
$\psi=\xi^{-1}\phi\xi$ which is isomorphic to $\phi$. By the two last cases and \cite[Proposition 2]{Po} we conclude
that $\hat{h}_{v,\phi}(\alpha)=\hat{h}_{v,\psi}(\alpha)\ge\frac1d$.

By the non-negativity of the local canonical heights we conclude that $\hat{h}_{\phi}(\alpha)\ge\frac1d$.
\end{proof}

\subsection{The general case}

The Lehmer problem for Drinfeld modules can be formulated in a more general set-up and its proof is reduced to that of
Theorem \ref{thmdm}. Let $C$ be a smooth projective irreducible curve defined over a finite field $\bbf_q$ of $q$
elements. Let $\infty$ be a fixed place of $K=\bbf_q(C)$, $A$ the ring of functions in $K$ which are regular everywhere
except at $\infty$, $v_{\infty}:K\to\mathbb{Z}\cup\{\infty\}$ the normalized discrete valuation associated to $\infty$
and $d_{\infty}$ the degree of $\infty$. For any $a\in A$, let $\deg(a)=-d_{\infty}v_{\infty}(a)$. The field $K$ is an
$A$-module with respect to the inclusion $A\subset K$. A Drinfeld $A$-module of rank $r$ defined over $K$ is a ring
homomorphism $\phi:A\to\text{End}_K(\mathbb{G}_a)\cong K\{\tau\}$ such that for every $a\in A$,
$\deg(\phi_a)=q^{r\deg(a)}$ and the constant term of $\phi_a$ is $a$ itself.

Let $a\in A-\bbf_q$. The global height of $\alpha\in\ov{K}$ is defined as (cf. \cite[\S2]{De})
$$\ch_{\phi}(\alpha)=\lim_{n\to\infty}\frac{h(\phi_{a^n}(\alpha))}{\deg(\phi_{a^n})}.$$
Let $L=K(\alpha)$ and $d=[L:K]$. For every discrete valuation $v:L\to\bbz\cup\{\infty\}$ of degree $d_v$ the local
height is defined as (cf. \cite[\S4]{Po})
$$\ch_{\phi,v}(\alpha)=\lim_{n\to\infty}-\frac{d_{v}}{d}
\frac{\min\{0,v(\phi_{a^n}(\alpha))\}}{\deg(\phi_{a^n})}.$$ As observed in \cite[Proposition 3]{Po} these heights are
independent of the choice of $a\in A-\bbf_q$.

The Dedekind domain $A$ is a finitely generated $\bbf_q$-algebra. Let $\caa$ the the set of generators of $A$ as an
$\bbf_q$-algebra, $T\in\caa$, $\deg(T)=d_T$ and $\phi_T=T+a_1(T)\tau+\ldots+a_{rd_T}(T)\tau^{rd_T}$. Let
$\alpha\in\ov{K}$ be a non-torsion element for $\phi$ with $[K(\alpha):K]=d$. Replacing the $a_i$'s in the proof of
Theorem \ref{thmdm} by the $a_i(T)$'s the proof of Theorem \ref{thmdm} shows that
\begin{equation}\label{amod}
\ch_{\phi}(\alpha)\ge\frac1d.
\end{equation}

\section{Elliptic curves}
Let $C$ be a smooth irreducible projective curve defined over an algebraically closed field $l$ of characteristic 0 or
$p>3$, let $K=l(C)$ be its function field and $\ov{K}$ its algebraic closure. Let $E/K$ be a non-constant semistable
elliptic curve defined over $K$, $\varphi_{\cae}:\cae\to C$ its minimal semi-stable regular model,
$j_{\cae}:C\to\bbp^1$ the $j$-map induced by $\varphi_{\cae}$ and $\ch_E:E(\ov{K})\to\bbr$ its canonical height.

Let $P\in E(\ov{K})$ and $L=K(P)$ the field generated by $K$ and the coordinates of $P$. Let $d=[L:K]$, $M_L$ the set
of places $v$ of $L$ which are normalized so that $v:L\to\mathbb{Z}\cup\{\infty\}$ is the corresponding discrete
valuation. Let $L_v$ be the completion of $L$ with respect to $v$ and $\lambda_v:E(K_v)\to\bbr$ its local N\'eron
function \cite[Chapter VI]{Si2}. Let $w=v_{|K}$, $e(v|w)$ the ramification index of $v$ over $w$,
$w'=e(v|w)w:K\to\mathbb{Z}\cup\{\infty\}$ the normalization of $w$, $K_w$ the completion of $K$ with respect to $w$ and
$n(v|w)=[L_v:K_w]$. Then
\begin{equation}\label{locglob1}
\ch_E(P)=\frac 1d\sum_{v\in M_L}n(v|w)\lambda_v(P),
\end{equation}
\cite[VI, Theorem 2.1]{Si2}.

Let $\frD_{E/K}$ be the minimal discrimant of $E/K$ and $d_{E/K}=\deg(\frD_{E/K})$. Since $E/K$ is semi-stable, it
follows from \cite[Chapter VII, Proposition 5.1]{Si1} that  \linebreak $w'(\frD_{E/K})=-w'(j_{\cae})$ for every pole
$w'$ of $j_{\cae}$, thus $\deg(j_{\cae})=d_{E/K}$. For every $v\in M_L$, let $v^+=\max\{v,0\}$.

\begin{lemma}\label{lem1}
\cite[Proposition 1.3]{HiSi} Let $A,N\ge 1$ be integers, $Q_0,\cdots,Q_{6AN}\in E(L_{v})$ distinct points. Then there
exists $P_0\cdots,P_N\in\{Q_0,\cdots,Q_{6AN}\}$ such that for each $i\ne l$,
$$\lambda_{v}(P_i-P_l)\ge\frac{1-A^{-1}}{12}{v}^+(j_{\cae}^{-1}).$$
\end{lemma}

\begin{proposition}\label{prop1}
$\#\{Q\in E(L)\,|\,\ch_E(Q)<\frac{d_{E/K}}{96d}\}\le 24$.
\end{proposition}

\begin{proof}
Denote $\cas=\{Q\in E(L)\,|\,\ch_E(Q)<\frac{d_{E/K}}{96d}\}$ and suppose $\#\cas>24$. Let $A=2$ and
$N+1=[\frac{\#\cas}{12}]>1$ the integral part of $\frac{\#\cas}{12}$, then $1<N+1\le\frac{\#\cas}{12}$. So we can
choose $12N+1$ distinct points $Q_0,\cdots,Q_{12N+1}$ in $\cas$. By Lemma \ref{lem1} there exist
$P_0,\cdots,P_N\in\{Q_0,\cdots,Q_{12N+1}\}$ such that $\lambda_v(P_i-P_l)\ge\frac1{24}v^+(j_{\cae}^{-1})$ for $i\ne l$.
It follows from the triangle inequality that
\begin{equation}\label{ht1}
H=\max_{Q\in\cas}\ch_E(Q)\ge\max_{1\le i\le N}\ch(P_i)\ge\frac 1{4N(N+1)}\sum_{i\ne l}\ch_E(P_i-P_l).
\end{equation}
Hence, by (\ref{locglob1}) and (\ref{ht1}),
\begin{equation}\label{ht2}
\begin{aligned}
H&\ge\frac 1{4N(N+1)d}\sum_{i\ne l}\sum_{v\in M_L}n(v|w)\lambda_v(P_i-P_l)\ge\frac1{96d}\sum_{v\in
M_L}n(v|w)v^+(j_{\cae}^{-1})\\
&=\frac1{96d}\sum_{w'\in M_K}\sum_{v|w}\frac{n(v|w)}{e(v|w)}{w'}^+(j_{\cae}^{-1})\ge\frac1{96d}\sum_{w'\in
M_K}{w'}^+(j_{\cae}^{-1})=\frac{d_{E/K}}{96d}.
\end{aligned}
\end{equation}
\end{proof}

\begin{remark}
We used the fact that $l$ is algebraically closed just to ensure that the poles of $j_{\cae}$ have all degree 1.
\end{remark}

As a consequence of Proposition \ref{prop1} we obtain a theorem which simultaneously deals with the Lehmer and the Lang
problems for elliptic curves over function fields. Recall that the Lang problem is to find a constant $c>0$ depending
on $E/K$ such that for every non-torsion point $P\in E(K)$ we have $\ch_E(P)\ge cd_{E/K}$.

\begin{theorem}\label{thmec}
Let $P\in E(\ov{K})$ be a non-torsion point of $E/K$ and $d=[K(P):K]$. Then there exists an absolute real constant
$c>0$ such that $\ch_E(P)\ge c\frac{d_{E/K}}d$.
\end{theorem}

\begin{proof}
Suppose $\ch_E(P)<\frac {d_{E/K}}{60000d}$. Then for every $1\le n\le25$, $\ch_E(nP)=n^2\ch_E(P)<\frac{d_{E/K}}{96d}$,
which contradicts Proposition \ref{prop1}. So we take $c=\frac1{60000}$.
\end{proof}

\begin{remark}
The constant for the Lehmer problem is $\frac{d_{E/K}}{60000}$ so it depends only on $\deg(j_{\cae})=d_{E/K}$, in the
semi-stable case.
\end{remark}

\begin{remark}\label{rem1}
In \cite[Theorem 0.2]{HiSi1} Hindry and Silverman proved Lang's conjecture for function fields over algebraically
closed fields of characteristic 0. In the case where $d_{E/K}\ge 24(g-1)$, where $g$ denotes the genus of $K$, they
obtained an absolute constant $c$. However, our constant is greater than theirs, thus improving the result. In the case
where $d_{E/K}<24(g-1)$, their constant depends exponentially on $g$, whereas ours is absolute and improves the
constant part of their bound. Nevertheless, we have just proved Lang's conjecture in the case of semi-stable elliptic
curves. Inspired on \cite[Theorem 0.2]{HiSi1} we had previously proved Lang's conjecture for semi-stable elliptic
curves over function fields of positive characteristic \cite[Theorem 5]{Pa} using \cite[Proposition 1.2]{HiSi}. First,
the bounds we obtained there do not have absolute constants, they depended not only on $g$ but also on the inseparable
degree of $j_{\cae}$. Furthermore, the present bound improves their constant parts. The reason for obtaining an
absolute constant is that \cite[Proposition 1.3]{HiSi} gives a lower bound which depends only on the choice of a
positive integer $A$, however the lower bound of \cite[Proposition 1.2]{HiSi} depends on the number $N+1$ of points
$P_0,\cdots,P_N$ chosen in $E(L_v)$ (cf. \cite[proof of Proposition 3]{Pa}).
\end{remark}

Another consequence of Proposition \ref{prop1} is a bound for the order of the torsion group $E(K)_{\tor}$.

\begin{corollary}\label{cor1}
$\#E(K)_{\tor}\le 24$.
\end{corollary}

\begin{remark}
Previous bounds for the torsion of elliptic curves over function fiels in characteristic 0 were obtained in
\cite[Theorem 7.2]{HiSi1} and in the case of characteristic $p$, Goldfeld and Szpiro treated the case where $C$ is
defined over a finite field \cite[Theorem 13]{GoSz}, but the result extends to algebraically closed fields of
characteristic $p>3$ and we also obtained a bound (cf. \cite[Theorem 7]{Pa}) using \cite[Proposition 3]{Pa}. In the
case of characteristic 0, the upper bound depended on $d_{E/K}$. Using Szpiro's theorem on the minimal discriminant of
elliptic curves over function fields \cite[Th\'eor\`eme 1]{Sz}, i.e., $d_{E/K}\le 6p^e(2g-2+f_{E/K})$, where $p^e$ is
the inseparable degree of $j_{\cae}$ and $f_{E/K}$ is the degree of the conductor divisor of $E/K$, it follows an upper
bound whose constant part is worse than the bound of Corollary \ref{cor1}. The bounds in characteristic $p$ (in the
semi-stable case) were $\sigma_{E/K}^2$, respectively $2\sigma_{E/K}^2$, where $\sigma_{E/K}=\frac{d_{E/K}}{f_{E/K}}$.
If $d_{E/K}\ge24p^e(g-1)$, then (using again \cite[Th\'eor\`eme 1]{Sz}) $\sigma_{E/K}\le12p^e$ and otherwise
$\sigma_{E/K}\le d_{E/K}<24p^e(g-1)$. Not only is the bound of Corollary \ref{cor1} absolute, but also it is better
than the estimates for $\sigma_{E/K}^2$.
\end{remark}

\subsection{Integral points}

Theorem \ref{thmec} and Corollary \ref{cor1} imply as in \cite[\S8]{HiSi1} an upper bound for the number of integral
points of an $S$-minmal Weierstrass equation of $E/K$.

Let $S$ be a finite set of places of $K$ and $R_S\subset K$ the ring of $S$-integers. For every $a\in K$ let
$h_K(a)=[K:l(a)]$. A Weierstrass equation $y^2=x^3+Bx+C$ with discriminant $\Delta$ is called $S$-minimal if
$h_K(\Delta)$ is minimal subject to $f(x)\in R_S[x]$. Let $\delta=\min\{\ch_E(P)\,|\,P\in E(K)-(E(K)_{\tor}\cap
E(R_S))\}$ and $\epsilon=\max\{\ch_E(P)\,|\,P\in E(R_S)\}$. In \cite[Lemma 1.2 (a)]{Si3} it is shown that
$\#E(R_S)\le\#E(K)_{\tor}(K)(1+2\sqrt{\frac{\epsilon}{\delta}})^{r_E}$, where $r_E$ denotes the rank of $E(K)$. It
follows from \cite[Remark 14]{Pa} that
\begin{equation}\label{bd2}
\epsilon\le p^e(12g+4\#S+5d_{E/K}).
\end{equation}

\begin{theorem}\label{thmintpt}
Let $y^2=x^3+Bx+C$ be an $S$-minimal Weierstrass equation for $E/K$. If $d_{E/K}\ge24p^e(g-1)$, then
$\#E(R_S)\le24(2299\sqrt{p^e\#S})^{r_E}$, otherwise $\#E(R_S)\le24(2021\sqrt{gp^e\#S})^{r_E}$.
\end{theorem}

\begin{proof}
By Theorem \ref{thmec}, $\delta\ge\frac{d_{E/K}}{60000}$. If $d_{E/K}\ge24p^e(g-1)$, then
$g\le\frac{d_{E/K}}{24p^e}+1$. Thus, since $\#S\ge 1$,
\begin{equation}\label{bd1}\begin{aligned}
\frac{\epsilon}{\delta}&\le60000\frac{p^e}{d_{E/K}}(12g+4\#S+5d_{E/K})
\\&\le60000p^e\left(12\left(\frac1{24p^e}+1\right)+9\#S\right)\\&\le1320000p^e\#S.
\end{aligned}\end{equation}
The first statement follows from (\ref{bd1}) and Corollary \ref{cor1}.

Suppose now that $d_{E/K}<24p^e(g-1)$. In this case, since $\#S\ge 1$ and $g\ge 2$, we have
\begin{equation}\label{bd3}\begin{aligned}
\frac{\epsilon}{\delta}&\le60000\frac{p^e}{d_{E/K}}(12g+4\#S+5d_{E/K})\le60000p^e(12g+9\#S)\\&\le1020000p^eg\#S.
\end{aligned}\end{equation}
The second statement follows from (\ref{bd3}) and Corollary \ref{cor1}.
\end{proof}

\begin{remark}
The bound of Theorem \ref{thmintpt} improves the bounds of \cite[Theorem 8.1]{HiSi1} in the case of characteristic 0
when $d_{E/K}\ge24(g-1)$. When $d<24(g-1)$ we also have an improvement of the constant part (which does not depend on
the rank $r_E$ of $E(K)$) if $g\ge3$. Note that in this latter case, the constant part of their bound depends on $g$,
whereas ours does not. In both cases the bound of Theorem \ref{thmintpt} improves that of \cite[Theorem 15]{Pa}.
\end{remark}

\subsection{Lehmer problem : the general case}
If we no longer suppose that $E/K$ is a semi-stable elliptic curve, then $\deg(j_{\cae})<d_{E/K}$ (cf. \cite[Chapter
VII, Proposition 5.1]{Si1}). In this case, instead of Proposition \ref{prop1}, we need to bound the cardinality of a
smaller set
\begin{equation}\label{bd4}
\#\left\{Q\in E(L)\,;\,\ch_E(Q)<\frac{\deg(j_{\cae})}{96d}\right\}\le 24.
\end{equation}
As a consequence, Theorem \ref{thmec} is replaced by: for every non-torsion point $P$ of $E$ of degree $d$ over $K$ we
have $\ch_E(P)\ge\frac{c'}d$, where $c'=\frac{\deg(j_{\cae})}{60000}$. We cannot obtain Lang's conjecture as in Theorem
\ref{thmec}, because Lemma \ref{lem1} involves $v^+(j_{\cae}^{-1})$, hence the proof of Proposition \ref{prop1} only
gives $\deg(j_{\cae})$ and not $d_{E/K}$.

\subsection{Integral points : the general case}
The bound of (\ref{bd4}) also implies that $\#E(K)_{\tor}\le24$. Note that in the general case
$$f_{E/K}<2\#\{\text{poles of }j_{\cae}\}\le2\deg_s(j_{\cae}),$$
where $\deg_s(j_{\cae})$ denotes the separable degree of $j_{\cae}$.

\begin{theorem}\label{intptgen}
Let $y^2=x^3+Bx+C$ be an $S$-minimal Weierstrass equation for $E/K$. If $d_{E/K}\ge24p^e(g-1)$, then
$\#E(R_S)\le24(13788\sqrt{gp^e\#S})^{r_E}$, otherwise $\#E(R_S)\le24(12121g\sqrt{p^e\#S})^{r_E}$.
\end{theorem}

\begin{proof}
If $d_{E/K}\ge24p^e(g-1)$, then
\begin{equation}\label{bd5}\begin{aligned}
\frac{\epsilon}{\delta}&\le60000\frac{p^e}{\deg(j_{\cae})}(12g+4\#S+5d_{E/K})\\
&\le60000\frac{p^e}{\deg(j_{\cae})}
\left(\frac{d_{E/K}}{2p^e}+12+9d_{E/K}\#S\right)\\
&\le1320000\frac{p^e}{\deg(j_{\cae})}d_{E/K}\#S.
\end{aligned}\end{equation}
Szpiro's discriminant theorem \cite[Th\'eor\`eme 1]{Sz} was first proved in the case of semi-stable elliptic curves.
However, this result was extended by Pesenti and Szpiro to any elliptic curve \cite[Th\'eor\`eme 0.1]{PeSz}. It follows
from \cite[Th\'eor\`eme 0.1]{PeSz}, (\ref{bd5}) and $f_{E/K}<2\deg_s(j_{\cae})$ that
\begin{equation}\label{bd6}\begin{aligned}
\frac{\epsilon}{\delta}&\le7920000\frac{p^{2e}}{\deg(j_{\cae})}(2g-2+f_{E/K})\#S\\
&\le23760000\frac{p^{2e}}{\deg(j_{\cae})}gf_{E/K}\#S\le47520000p^eg\#S.
\end{aligned}\end{equation}
The result now follows from (\ref{bd6}) and $\#E(K)_{\tor}\le24$.

Suppose now that $d_{E/K}<24p^e(g-1)$, then
\begin{equation}\label{bd7}\begin{aligned}
\frac{\epsilon}{\delta}&\le60000\frac{p^e}{\deg(j_{\cae})}(12g+4\#S+5d_{E/K})
\le1020000\frac{p^e}{\deg(j_{\cae})}gd_{E/K}\#S\\
&\le6120000\frac{p^{2e}}{\deg(j_{\cae})}g(2g-2+f_{E/K})\#S\le18360000\frac{p^{2e}}{\deg(j_{\cae})}g^2f_{E/K}\#S\\
&\le36720000p^eg^2\#S.
\end{aligned}\end{equation}
The result now follows from (\ref{bd7}) and $\#E(K)_{\tor}\le24$.
\end{proof}

\begin{remark}
Observe that the bounds of Theorem \ref{intptgen} are a worse than those of Theorem \ref{thmintpt}.
\end{remark}

\end{document}